# QUASIFOLDS, DIFFEOLOGY AND NONCOMMUTATIVE GEOMETRY

PATRICK IGLESIAS-ZEMMOUR AND ELISA PRATO

ABSTRACT. After embedding the objects quasifolds into the category {Diffeology}, we associate a $\mathbf{C}^*$-agebra with every atlas of any quasifold, and show how different atlases give Morita equivalent algebras. This builds a new bridge between diffeology and noncommutative geometry — beginning with the today classical example of the irrational torus — which associates a Morita class of $\mathbf{C}^*$-algebras with a diffeomorphic class of quasifolds.

## Introduction

This paper is a follow-up to "Noncommutative geometry & diffeology: the case of orbifolds" [IZL17]. In that article, a construction was established that associated a $\mathbf{C}^*$-algebra with every orbifold, in a functorial way. Here we extend the construction to the more general *quasifolds*.[1]

First of all, we identify quasifolds [EP01] as objects in the category {Diffeology} [PIZ13]. These are *diffeological spaces* that are locally diffeomorphic, at each point, to some quotient $\mathbf{R}^n/\Gamma$, for some integer $n$, and for $\Gamma$ — which may change from point to point — a countable subgroup of $\mathrm{Aff}(\mathbf{R}^n)$. As it appears clearly, the definition is similar to that of orbifolds [IKZ10], except for the group $\Gamma$, which can be infinite, while it is finite for orbifolds — whose original definition as independent objects has been published by Ishiro Satake in [IS56, IS57].

By considering quasifolds as diffeologies, they inherit a structure of category which we denote by {Quasifolds}, whose morphisms are smooth maps in the sense of diffeology. This remark carries a strong content, as can be seen in the lifting of smooth functions between quasifolds (§2): the same phenomenon happening for orbifolds, where smooth maps may not lift locally equivariantly, happens also for strict[2] quasifolds. We have *a*

---

*Date*: August 3, 2020.

2010 *Mathematics Subject Classification*. Primary 53; Secondary 58B34.

Research supported by the Laboratory Ypatia of Mathematical Sciences LIA LYSM AMU-CNRS-ECM-INdAM. And for E.P. also supported by the PRIN Project "Real and Complex Manifolds: Geometry, Topology and Harmonic Analysis" (MIUR, Italy) and by GNSAGA (INdAM, Italy).

[1]This also provides an answer to the referee of the first paper, who suggested to move on to more general diffeologies, in order to produce more interesting objects.

[2]Which are not orbifolds.

*priori* the sequence of categories:

$$\{\text{Manifolds}\} \preceq \{\text{Orbifolds}\} \preceq \{\text{Quasifolds}\} \preceq \{\text{Diffeology}\}.$$

Then, we generalize to quasifolds the functor toward noncommutative geometry, developed in [IZL17] for orbifolds. In the same way, we associate with each atlas of a quasifold a *structure groupoid*,[3] in (§5). The objects of this groupoid are the elements of the nebula of the strict generating family associated with the atlas. The arrows between the objects are the germs of the local diffeomorphisms of the nebula that are absorbed by the evaluation map, that is, which project to the identity on the quasifold.

In parallel with the case of orbifolds, in (§3) and in (§4) we generalize to quasifolds the two fundamental results:

THEOREM. *Any local smooth map on $\mathbf{R}^n$ that projects to the identity in the quotient $\mathbf{R}^n/\Gamma$, where $\Gamma$ is a countable subgroup of $\mathrm{Aff}(\mathbf{R}^n)$, is everywhere locally the action of some $\gamma \in \Gamma$.*

THEOREM. *Local diffeomorphisms between quasifolds lift by local diffeomorphisms on the level of the strict generating families. Pointed local diffeomorphisms lift by pointed local diffeomorphisms, where the source and the target can be chosen arbitrarily in the appropriate fibers over the quasifold.*

The difficulty here is to pass from the action of a finite group on a Euclidean domain to the action of a possibly infinite, but countable group, whose orbits can be dense. This has led us to a substantial revision of the methods, focusing on the countable nature of the groups, and has resulted in proofs that are minimal and essential.

As said above, in (§5) we define the structure groupoid associated with an atlas of the quasifold. Then, thanks to the previous theorem, in (§6) we prove the following:

THEOREM. *Two different atlases of a same quasifold give two equivalent groupoids, as categories* [SML78]. *Consequently, two diffeomorphic quasifolds have equivalent structure groupoids.*

In other words, the class of the structure groupoid is a diffeological invariant of the quasifold. Then, in (§7), we give a general description of the structure groupoids.

Next, in (§9) we prove the following:

THEOREM. *The groupoids associated with two different atlases of a same quasifold satisfy the Muhly-Renault-Williams equivalence.*

Then, having proved in (§8) that the structure groupoids associated with the atlases of a quasifold are étale and Hausdorff, we show that they fulfill the conditions of Jean Renault's construction of an associated $\mathbf{C}^*$-algebra, by equipping the set of morphisms with the same counting measure as in the case of orbifolds. And in (§10) we prove then, thanks to (§9), the main result:

---

[3] Which can be regarded as the *gauge groupoid* of the quasifold structure.

THEOREM. *The $\mathbf{C}^*$-algebras associated with different atlases of a same quasifold are Morita-equivalent. Therefore, diffeomorphic quasifolds have Morita-equivalent $\mathbf{C}^*$-algebras.*

Finally, we illustrate this construction with two simple examples: the traditional irrational torus $T_\alpha$ and the $\mathbf{Q}$-circle, quotient of $\mathbf{R}$ by $\mathbf{Q}$. In these two examples, we observe that our construction gives the expected result. In work in progress, we apply these techniques also to the class of *symplectic toric quasifolds* [EP01, FBEP01].

From the very beginning, with the 1983 paper [PDPI83] on the irrational torus $T_\alpha$, it was clear that there existed some connection between diffeology and noncommutative geometry. Beginning with the fact that two such tori $T_\alpha$ and $T_\beta$ were diffeomorphic if and only if $\alpha$ and $\beta$ were equivalent modulo $\mathrm{GL}(2,\mathbf{Z})$, which is the same condition for their algebra to be Morita-equivalent [MR81]. That could not be just chance. This work, which began with the case of orbifolds [IZL17] and which continues here with quasifolds, shows and describes the logic behind this correspondence. We can reasonably expect wider links between the two theories, which will be addressed in the future.

NOTE 1. Unlike the categorical approach, which defines its objects directly by means of higher structures (stacks, n-categories etc.), we induce the groupoid generating the $\mathbf{C}^*$-algebra of the quasifold via its singular geometry encoded in the diffeology. So, to the current standard way {groupoid → $\mathbf{C}^*$-algebra}, we add a first floor {diffeology → groupoid}, which is not trivial and makes this construction non-tautological.

NOTE 2. The irrational tori in arbitary dimension, or quasitori, are particular quasifolds that are dual smooth geometric versions of quasilattices. Knowing how to associate a $\mathbf{C}^*$-algebra to a quasifold in a structural fashion can be viewed as a kind of geometric quantization. In fact, the study of the spectrum of the Hamiltonian in a quasicrystal was at the origin of Alain Connes' noncommutative geometry. It is obviously interesting to have a smooth version of this, which is what we are providing.

NOTE 3. We assume that the reader is familiar with the basic concepts in diffeology and we refer to the textbook [PIZ13] for details. Let us just recall that a diffeology on a set X is a set $\mathscr{D}$ of smooth parametrizations, called *plots*, that satisfy three fundamental axioms: *covering*, *locality* and *smooth compatibility*. That said, there are a couple of important diffeological constructions that we use in the following. First, the *quotient diffeology*: every quotient of a diffeological space inherits a natural diffeology for which the plots are the parametrizations that can be locally lifted by plots in the source space. Then, the *subset diffeology*: every subset of a diffeological space inherits a diffeology for which the plots are the plots of the ambient space, but with values in the subset. For example, in diffeology a subset is *discrete* if the subset diffeology is the discrete diffeology, that is, the plots are locally constant. For example $\mathbf{Q} \subset \mathbf{R}$ is discrete. Finally, the *local diffeology*:[4] a map $f$, from a subset A of a diffeological space X to a diffeological space $\mathrm{X}'$, is a *local smooth map* if and only if its composite $f \circ \mathrm{P}$ with a plot P in X, defined on $\mathrm{P}^{-1}(\mathrm{A})$, is a plot in

---

[4]Introduced with the definition of local smooth maps in [PI85, §1.2.3], see also [PIZ13, §2.1].

X′. That is equivalent to: A is an open subset for the *D-topology*[5] and $f$ restricted to A is smooth for the subset diffeology. With local smooth maps come local diffeomorphisms, which are the fundamentals of modeling spaces in diffeology [PIZ13, §4.19], on which many constructions of subcategories are based, like manifolds, manifolds with boundary and corners, orbifolds, quasifolds etc.

THANKS. It is a pleasure for one of the authors (PIZ) to thank Anatole Khelif for useful discussions on $\mathbf{C}^*$-algebras.

## Diffeological Quasifolds

The notion of *quasifold* has been introduced in 1999 in the paper "On a generalization of the notion of orbifold" [EP99], see also [EP01]. The idea is that a *n*-quasifold is a smooth object which resembles locally everywhere a quotient $\mathbf{R}^n/\Gamma$, where $\Gamma$ is some countable subgroup of diffeomorphisms. The analogy with orbifolds, for which $\Gamma$ is finite, is indeed clear. On the other hand, *Diffeology* has been precisely developed, from the mid '80, to deal with this kind of situation, beginning with "Exemple de groupes différentiels…" [PDPI83]. In particular, *orbifolds* have been later successfully included as a subcategory in {Diffeology} in the paper "Orbifolds as Diffeology" [IKZ10]. It was natural to try to include also quasifolds, and this is what we do now.

1. WHAT IS A DIFFEOLOGICAL QUASIFOLD? — We have indeed a diffeological version of quasifolds, formally defined by:

DEFINITION. *A n-quasifold is a diffeological space* X *which is locally diffeomorphic, everywhere, to some* $\mathbf{R}^n/\Gamma$*, where* $\Gamma$ *is a countable subgroup, maybe infinite, of* $\mathrm{Aff}(\mathbf{R}^n)$. *The group* $\Gamma$ *maybe changing from place to place.*

In more words, this definition means precisely the following: for all $x \in X$, there exist a countable subgroup $\Gamma \subset \mathrm{Aff}(\mathbf{R}^n)$, and a local diffeomorphism $\varphi$ from $\mathbf{R}^n/\Gamma$ to X, defined on some open subset $U \subset \mathbf{R}^n/\Gamma$, such that $x \in \varphi(U)$. The subset U is open for the D-topology, that is in this case, the quotient topology [PIZ13, §2.12] by the projection map[6] class: $\mathbf{R}^n \to \mathbf{R}^n/\Gamma$. That said:

DEFINITION. *Any such diffeomorphism is called a* chart. *A set of charts* $\mathscr{A}$*, covering* X*, is called an* atlas.

NOTE. In the following we consider only quasifolds that support a *locally finite* atlas, that is, every point in the quasifold is covered by a finite number of charts. For example, a symplectic toric quasifold has a canonical atlas made of a finite number of charts [FBEP01, Thm. 3.2].

REMARK 1. This approach to quasifolds considers spaces that are already equipped with a smooth structure, that is, a diffeology, and then, checks if that diffeology is generated by local diffeomorphisms with some quotients $\mathbf{R}^n/\Gamma$. This is the standard construction of modeling diffeology we mentioned above; it applies to manifolds, orbifolds…and now quasifolds. It is a reverse construction as the usual one, where the smooth structure is built after equipping the underlying set with a family of injections, compatible according to some specific conditions. Recent works and results involving quasifolds in symplectic geometry can be found in [FBEP18], [FBEP19], and [BPZ19].

---

[5]See [PIZ13, §2.8].

[6]In this paper the word class will denote generically the *class map* from a space onto its quotient, for a relation which has been clearly identified.

REMARK 2. The group $\Gamma$ is chosen inside the affine group and not just the linear subgroup, as it is the case for orbifolds. In this way, one immediately has the well known example of the *irrational torus* $T_\alpha = \mathbf{R}/\mathbf{Z} + \alpha\mathbf{Z}$ [PDPI83], where $\alpha x \in \mathbf{R} - \mathbf{Q}$, as a quasifold. But, we can notice that $\Gamma$ could be embedded in $\mathrm{GL}(n+1,\mathbf{R})$ by considering $\mathbf{R}^n$ as the subspace of height 1 in $\mathbf{R}^n \times \mathbf{R}$, and an element $(A,b) \in \mathrm{Aff}(\mathbf{R}^n)$ acting on $\mathbf{R}^n \times \{1\}$ by

$$\begin{pmatrix} A & b \\ 0 & 1 \end{pmatrix} \begin{pmatrix} X \\ 1 \end{pmatrix} = \begin{pmatrix} AX + b \\ 1 \end{pmatrix}.$$

Hence, the affine or linear nature for the subgroup $\Gamma$ is not really discriminant.

REMARK 3. In *Example of Singular Reduction in Symplectic Diffeology* [PIZ16], an infinite dimensional quasi-projective space is built inside the category of diffeology. That is, an example of an infinite dimensional analog of the present concept of quasifold. That leaves some space for a generalization of the kind of constructions explored in this paper.

2. SMOOTH MAPS BETWEEN QUASIFOLDS — As an object of the category of diffeological spaces, quasifolds inherit automatically the notion of smooth maps. A smooth map from a quasifold to another quasifold is just a map which is smooth when the quasifolds are regarded as diffeological spaces. It follows immediately that the composite of smooth maps between quasifolds is again a smooth map. Hence, quasifolds form a full subcategory of {Diffeology} we shall denote by {Quasifolds}.

A special phenomenon appearing in the case of orbifolds persists for quasifolds: smooth maps between diffeological quasifolds may have no local equivariant lifting, as shown by the following example inspired by [IKZ10, Example 25].

Let $\alpha \in \mathbf{R} - \mathbf{Q}$ and $\mathbf{C}_\alpha$ be the irrational quotient:

$$\mathbf{C}_\alpha = \mathbf{C}/\Gamma \quad \text{with} \quad \Gamma = \{e^{i2\pi\alpha k}\}_{k \in \mathbf{Z}}.$$

This diffeological space[7] falls into the category of quasifolds.

Let now $f : \mathbf{C} \to \mathbf{C}$ be defined by

$$f(z) = \begin{cases} 0 & \text{if } r > 1 \text{ or } r = 0 \\ e^{-1/r} \varrho_n(r)\, r & \text{if } \frac{1}{n+1} < r \leq \frac{1}{n} \text{ and } n \text{ is even} \\ e^{-1/r} \varrho_n(r)\, z & \text{if } \frac{1}{n+1} < r \leq \frac{1}{n} \text{ and } n \text{ is odd,} \end{cases}$$

---

[7]Appearing already in [PI85, Appendix 6].

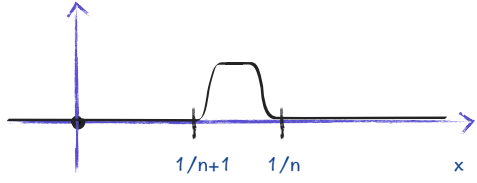

FIGURE 1. The function $\rho_n$.

where $r = \sqrt{|z|^2}$ and $\rho_n$ is a function vanishing flatly outside the interval $]1/(n+1), 1/n[$ and not inside, see Figure 1.

If we consider now $\tau \in U(1)$, one has: $f(\tau z) = f(z)$ on the annulus $\frac{1}{n+1} < r \leq \frac{1}{n}$ if $n$ even, and $f(\tau z) = \tau f(z)$ if $n$ is odd. That is, $f(\tau z) = h_z(\tau) f(z)$, where $h_z(\tau) = 1$ or $h_z(\tau) = \tau$ depending on whether $z$ is in an even or odd annulus. Hence, class$(f(\gamma z)) =$ class$(f(z))$ for all $\gamma \in \Gamma$. Then, the map $f$ projects onto a smooth map $\varphi \colon \mathbf{C}_\alpha \to \mathbf{C}_\alpha$ defined by

$$\varphi(\text{class}(z)) = \text{class}(f(z)).$$

Next, assume that $f'$ is another lifting of $\varphi$. For all $z \in \mathbf{C}$, there exists $\gamma(z) \in \Gamma$ such that $f'(z) = \gamma(z) f(z)$. We then get a smooth map $z \mapsto f'(z)/f(z) = \gamma(z)$ defined on $\mathbf{C} - \{0\}$. Since $\Gamma \subset U(1)$ is diffeologically discrete [PIZ13, Exercise 8 p. 14], this map is constant $\gamma(z) = \gamma$ and $f'(z) = \gamma f(z)$ on $\mathbf{C} - \{0\}$, and by continuity on $\mathbf{C}$. Thus, two lifts of $\varphi$ differ only by a constant in $\Gamma$, which gives the same function $h'_z = h_z$. Therefore, because the homomorphism $h_z$ flips from the trivial homomorphism to the identity on successive annuli, $\varphi$ has no local equivariant smooth lifting.

3. LIFTING THE IDENTITY — Let $\mathscr{Q} = \mathbf{R}^n/\Gamma$. Consider a local smooth map F from $\mathbf{R}^n$ to itself, such that class $\circ$ F = class. In other words, F is a local lifting of the identity on $\mathscr{Q}$. Then,

THEOREM. F *is locally equal to some group action* $F(r) =_{\mathrm{loc}} \gamma \cdot r = Ar + b$, *where* $\gamma = (A, b) \in \Gamma$, *for some* $A \in \mathrm{GL}(\mathbf{R}^n)$ *and* $b \in \mathbf{R}^n$.

*Proof.* Let us assume first that F is defined on an open ball $\mathscr{B}$. Then, for all $r$ in the ball, there exists a $\gamma \in \Gamma$ such that $F(r) = \gamma \cdot r$. Next, for every $\gamma \in \Gamma$, let

$$F_\gamma \colon \mathscr{B} \to \mathbf{R}^n \times \mathbf{R}^n \quad \text{with} \quad F_\gamma(r) = (F(r), \gamma \cdot r).$$

Let $\Delta \subset \mathbf{R}^n \times \mathbf{R}^n$ be the diagonal and let us consider

$$\Delta_\gamma = F_\gamma^{-1}(\Delta) = \{r \in \mathscr{B} \mid F(r) = \gamma \cdot r\}.$$

LEMMA 1. *There exist at least one* $\gamma \in \Gamma$ *such that the interior* $\mathring{\Delta}_\gamma$ *is non-empty.*

◂ Indeed, since $F_\gamma$ is smooth (thus continuous), the preimage $\Delta_\gamma$ by $F_\gamma$ of the diagonal is closed in $\mathscr{B}$. However, the union of all the preimages $F_\gamma^{-1}(\Delta)$ — when $\gamma$ runs over $\Gamma$ — is the ball $\mathscr{B}$. Then, $\mathscr{B}$ is a countable union of closed subsets. According to Baire's theorem, there is at least one $\gamma$ such that the interior $\mathring{\Delta}_\gamma$ is not empty. ▶

LEMMA 2. *The union $\mathring{\Delta}_\Gamma = \cup_{\gamma \in \Gamma} \mathring{\Delta}_\gamma$ is an open dense subset of $\mathscr{B}$.*

◀ Indeed, let $\mathscr{B}' \subset \mathscr{B}$ be an open ball. Let us denote with a prime the sets defined above but for $\mathscr{B}'$. Then, $\Delta'_\gamma = (F_\gamma \restriction \mathscr{B}')^{-1}(\Delta) = \Delta_\gamma \cap \mathscr{B}'$, and then $\mathring{\Delta}'_\gamma = \mathring{\Delta}_\gamma \cap \mathscr{B}'$. Thus, $\mathscr{B}' \cap \mathring{\Delta}_\Gamma = \mathscr{B}' \cap (\cup_{\gamma \in \Gamma} \mathring{\Delta}_\gamma) = \cup_{\gamma \in \Gamma} \mathring{\Delta}'_\gamma$, which is not empty for the same reason that $\cup_{\gamma \in \Gamma} \mathring{\Delta}_\gamma$ is not empty. Therefore, $\mathring{\Delta}_\Gamma$ is dense. ▶

Hence, there exists a subset of $\Gamma$, indexed by a family $\mathscr{I}$, for which $\mathcal{O}_i = \mathring{\Delta}_{\gamma_i} \subset \mathscr{B}$ is open and non-empty, $\cup_{i \in \mathscr{I}} \mathcal{O}_i$ is an open dense subset of $\mathscr{B}$, and $F \restriction \mathcal{O}_i : r \mapsto A_i r + b_i$, where $(A_i, b_i) \in \mathrm{Aff}(\mathbf{R}^n)$. Since F is smooth, the first derivative $D(F)$ restricted to $\mathcal{O}_i$ is equal to $A_i$, and then the second derivative $D^2(F) \restriction \mathcal{O}_i = 0$, for all $i \in \mathscr{I}$. Then, since $D^2(F) = 0$ on an open dense subset of $\mathscr{B}$, $D^2(F) = 0$ on $\mathscr{B}$, that is $D(F)(r) = A$ for all $r \in \mathscr{B}$, with $A \in \mathrm{GL}(n, \mathbf{R})$. Now, the map $r \mapsto F(r) - Ar$, defined on $\mathscr{B}$, is smooth. But, restricted on $\mathcal{O}_i$ it is equal to $b_i$. Its derivative vanishes on the open dense subset $\cup_{i \in \mathscr{I}} \mathcal{O}_i$ and thus vanishes on $\mathscr{B}$. Therefore, $F(r) - Ar = b$ on the whole $\mathscr{B}$, for $b \in \mathbf{R}^n$ and $F(r) = Ar + b$ on $\mathscr{B}$, with $\gamma = (A, b) \in \Gamma$. □

4. LIFTING LOCAL DIFFEOMORPHISMS — Let $\mathscr{Q} = \mathbf{R}^n / \Gamma$ and $\mathscr{Q}' = \mathbf{R}^{n'} / \Gamma'$, where $\Gamma \subset \mathrm{Aff}(\mathbf{R}^n)$ and $\Gamma' \subset \mathrm{Aff}(\mathbf{R}^{n'})$ are countable subgroups. Then,

THEOREM. *Every local smooth lifting $\tilde{f}$ of any local diffeomorphism $f$ of $\mathscr{Q}$ is necessarily a local diffeomorphism. In particular $n = n'$. Moreover, let $x \in \mathrm{dom}(f)$, $x' = f(x)$, $r, r' \in \mathbf{R}^n$ be such that $\mathrm{class}(r) = x$ and $\mathrm{class}(r') = x'$. Then, the local lifting $\tilde{f}$ can be chosen such that $\tilde{f}(r) = r'$.*

Note that $n$ is also the diffeological dimension of $\mathbf{R}^n / \Gamma$, see [PIZ13, §1.78].

*Proof.* Let the local diffeomorphism $f$ be defined on U with values in U$'$. By definition of local diffeomorphism, they are both open for the D-topology. Then $\tilde{U} = \mathrm{class}^{-1}(U)$ is open in $\mathbf{R}^n$. Since the composite $f \circ \mathrm{class} : \tilde{U} \to U'$ is a plot in $\mathscr{Q}'$, for all $r \in \tilde{U}$ there exists a smooth local lifting $\tilde{f} : \tilde{V} \to \mathbf{R}^{n'}$, defined on an open neighborhood of $r$, such that $\mathrm{class}' \circ \tilde{f} = f \circ \mathrm{class} \restriction \tilde{V}$.

$$\begin{array}{ccc} \mathbf{R}^n \supset \tilde{U} \supset \tilde{V} \xrightarrow{\tilde{f}} \mathbf{R}^{n'} & \quad & \mathbf{R}^n \xleftarrow{\hat{f}} \tilde{V}' \subset \tilde{U}' \subset \mathbf{R}^{n'} \\ \mathrm{class} \downarrow \qquad \qquad \downarrow \mathrm{class}' & \quad & \mathrm{class} \downarrow \qquad \qquad \downarrow \mathrm{class}' \\ \mathscr{Q} \supset U \xrightarrow{f} \mathscr{Q}' & \quad & \mathscr{Q} \xleftarrow{f^{-1}} U' \subset \mathscr{Q}' \end{array}$$

Let $x = \mathrm{class}(r)$, $x' = f(x)$, $r' = \tilde{f}(r)$, and then $x' = \mathrm{class}'(r')$.

Next, let $\tilde{U}' = \text{class}'^{-1}(U')$. Since the composite $f^{-1} \circ \text{class}'$ is a plot in $\mathcal{Q}$, there exists a smooth lifting $\hat{f}\colon \tilde{V}' \to \mathbf{R}^n$, defined on an open neighborhood of $r'$, such that $\text{class} \circ \hat{f} = f^{-1} \circ \text{class}' \restriction \tilde{V}'$. Let $r'' = \hat{f}(r')$, which is a priori different from $r$.

Now, we consider the composite $\hat{f} \circ \tilde{f}\colon \tilde{W} \to \mathbf{R}^n$, where $\tilde{W} = \tilde{f}^{-1}(\tilde{V}')$ is a non-empty open subset of $\mathbf{R}^n$ since it contains $r$. Moreover, $\hat{f} \circ \tilde{f}(r) = r''$. It also satisfies $\text{class} \circ (\hat{f} \circ \tilde{f}) = \text{class}$. Indeed, $\text{class} \circ (\hat{f} \circ \tilde{f}) = (\text{class} \circ \hat{f}) \circ \tilde{f} = (f^{-1} \circ \text{class}') \circ \tilde{f} = f^{-1} \circ (\text{class}' \circ \tilde{f}) = f^{-1} \circ (f \circ \text{class}) = (f^{-1} \circ f) \circ \text{class} = \text{class}$. Thus, thanks to (§3), there exists, locally, $\gamma \in \Gamma$ such that $\hat{f} \circ \tilde{f} = \gamma \restriction \tilde{W}$. By the way, $r'' = (\hat{f} \circ \tilde{f})(r) = \gamma \cdot r$. Let $\bar{f} = \gamma^{-1} \circ \hat{f}$, then: $\text{class} \circ \bar{f} = \text{class} \circ \gamma^{-1} \circ \hat{f} = \text{class} \circ \hat{f} = f^{-1} \circ \text{class}'$, and $\bar{f}$ is still a local lifting of $f^{-1}$. Thus $\bar{f} \circ \tilde{f} = \mathbf{1}_{\tilde{W}}$, that is, $\bar{f} = \tilde{f}^{-1} \restriction \tilde{W}$. We conclude that, around $r$, $\tilde{f}$ is a local diffeomorphism. Now, if we consider any another point $r'''$ over $x'$, there exists $\gamma'$ such that $\gamma' \cdot r' = r'''$; changing $\tilde{f}$ to $\gamma' \circ \tilde{f}$ and $\bar{f}$ to $\bar{f} \circ \gamma'^{-1}$, we get $\tilde{f}(r) = r'''$, and $\tilde{f}$ and $\bar{f}$ still remain inverse of each other.

Therefore, for any $r \in \mathbf{R}^n$ over $x$ and any $r' \in \mathbf{R}^n$ over $x' = f(x)$, we can locally lift $f$ to a local diffeomorphism $\tilde{f}$ such that $\tilde{f}(r) = r'$. $\square$

## Structure Groupoids of a Quasifold.

In this section, we associate a *structure groupoid* — or *gauge groupoid* — which is a diffeological groupoid [PIZ13, 8.3], with every atlas of a quasifold. Then we show that different atlases give equivalent groupoids: as categories, according to the Mac Lane definition [SML78], and in the sense of Muhly-Renault-Williams [MRW87]. We give a precise description of the structure groupoid in terms or the groupoid associated with the action of the structure groups $\Gamma$, and the connecting points of the charts. This construction is the foundation for a $\mathbf{C}^*$-algebra associated with the quasifold.

5. Building the groupoid of a quasifold. — Let X be a quasifold, let $\mathscr{A}$ be an atlas and let $\mathscr{F}$ be the strict generating family over $\mathscr{A}$. We denote by $\mathscr{N}$ the nebula[8] of $\mathscr{F}$, that is, the sum of the domains of its elements:

$$\mathscr{N} = \coprod_{F \in \mathscr{F}} \text{dom}(F) = \{(F, r) \mid F \in \mathscr{F} \text{ and } r \in \text{dom}(F)\}.$$

The *evaluation map* is the natural subduction

$$\text{ev}\colon \mathscr{N} \to X \quad \text{with} \quad \text{ev}(F, r) = F(r).$$

Following the construction in the case of orbifolds [IZL17], the *structure groupoid* of the quasifold X, associated with the atlas $\mathscr{A}$, is defined as the subgroupoid $\mathbf{G}$ of germs of local diffeomorphisms of $\mathcal{N}$ that project to the identity of X along ev. That is,

$$\begin{cases} \operatorname{Obj}(\mathbf{G}) &= \mathcal{N}, \\ \operatorname{Mor}(\mathbf{G}) &= \{\operatorname{germ}(\Phi)_\nu \mid \Phi \in \operatorname{Diff}_{\operatorname{loc}}(\mathcal{N}) \text{ and } \operatorname{ev} \circ \Phi = \operatorname{ev} \restriction \operatorname{dom}(\Phi)\}. \end{cases}$$

---

[8]See definition in [PIZ13, § 1.76].

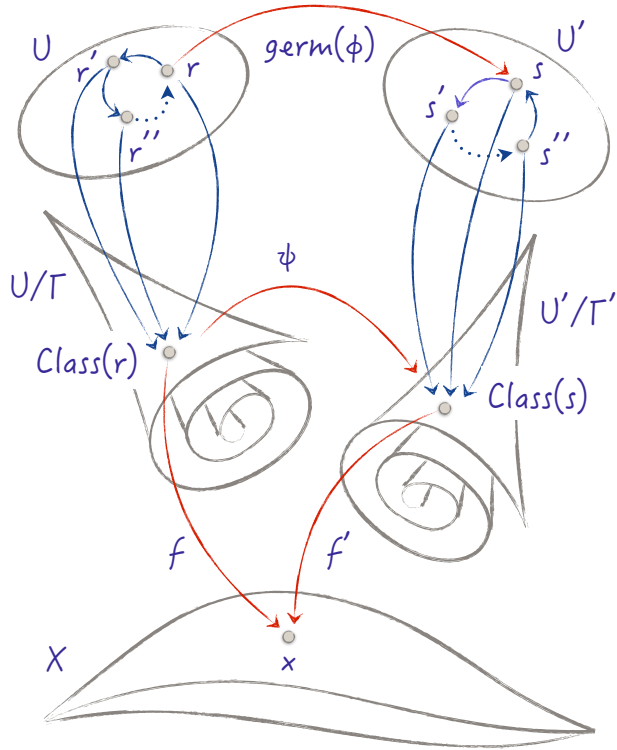

Figure 2. The three levels of a quasifold.

The set $\operatorname{Mor}(\mathbf{G})$ is equipped with the functional diffeology inherited by the full groupoid of germs of local diffeomorphisms [IZL17, §2 & 3]. Note that, given $\Phi \in \operatorname{Diff}_{\operatorname{loc}}(\mathcal{N})$ and $\nu \in \operatorname{dom}(\Phi)$, there exist always two plots F and F′ in $\mathcal{F}$ such that $\nu = (F, r)$, with $r \in \operatorname{dom}(F)$, and a local diffeomorphism $\varphi$ of $\mathbf{R}^n$, defined on an open ball centered in $r$, such that $\operatorname{dom}(\varphi) \subset \operatorname{dom}(F)$, $\varphi = \Phi \restriction \{F\} \times \operatorname{dom}(F)$ and $F' \circ \varphi = F \restriction \operatorname{dom}(\varphi)$. That is summarized by the diagram:

$$\operatorname{dom}(F) \supset \operatorname{dom}(\varphi) \xrightarrow{\varphi} \operatorname{dom}(F') \\ \phantom{\operatorname{dom}(F) \supset \operatorname{dom}(\varphi)} \searrow^{F} \quad \swarrow^{F'} \\ \phantom{\operatorname{dom}(F) \supset \operatorname{dom}(\varphi) \searrow^{F}} X$$

NOTE. According to the theorem in (§3), the local diffeomorphisms, defined on the domain of a generating plot, and lifting the identity of the quasifold, are just the elements of the structure group associated with the plot. We can legitimately wonder what is the

point of involving general germs of local diffeomorphisms, if we merely end up with the structure group we could have began with. The reason is that the structure groups connect the points of the nebula that project on a same point of the quasifold, only when they are inside the same domain. They cannot connect the points of the nebula that project on the same point of the quasifold but belonging to different domains, with maybe different structure groups. This is the reason why we cannot avoid the use of germs of local diffeomorphisms in the nebula, to begin with. That situation is illustrated in Figure 2.

6. EQUIVALENCE OF STRUCTURE GROUPOIDS Let us recall that a functor $S\colon A \to C$ is an equivalence of categories if and only if, S is full and faithful, and each object $c$ in C is isomorphic to $S(a)$ for some object $a$ in A [SML78, Chap. 4 § 4 Thm. 1]. If A and C are groupoids, the last condition means that, for each object $c$ of C, there exist an object $a$ of A and an arrow from $S(a)$ to $c$.

In other words: let the *transitivity-components* of a groupoid be the maximal full subgroupoids such that each object is connected to any other object by an arrow. The functor S is an equivalence of groupoids if it is full and faithful, and surjectively projected on the set of transitivity-components.

Now, consider an $n$-quasifold X. Let $\mathscr{A}$ be an atlas, let $\mathscr{F}$ be the associated strict generating family, let $\mathscr{N}$ be the nebula of $\mathscr{F}$ and let **G** the associated structure groupoid. Let us first describe the *morphology* of the groupoid.

PROPOSITION. *The fibers of the subduction* $\mathrm{ev}\colon \mathrm{Obj}(\mathbf{G}) \to \mathrm{X}$ *are exactly the transitivity-components of* **G**. *In other words, the space of transitivity components of the groupoid* **G** *associated with any atlas of the quasifold* X, *equipped with the quotient diffeology, is the quasifold itself.*

THEOREM. *Different atlases of* X *give equivalent structure groupoids. The structure groupoids associated with diffeomorphic quasifolds are equivalent.*

In other words, the equivalence class of the structure groupoids of a quasifold is a diffeological invariant.

*Proof.* These results are analogous to the results of [IZL17, §5]. They have the same kind of proof. The fact that the structure groups $\Gamma$ of the quasifolds are countable instead of finite has no negative consequences, thanks to (§4).

Let us start by proving the proposition. Let $F\colon U \to X$ and $F'\colon U' \to X'$ be two generating plots from the strict family $\mathscr{F}$, and $r \in U \subset \mathbf{R}$ and $r' \in U' \subset \mathbf{R}'$. Assume that $\mathrm{ev}(F, r) = \mathrm{ev}(F', r') = x$, that is, $x = F(r) = F'(r')$. Note that $F = f \circ \mathrm{class} \restriction U$ and $F' = f' \circ \mathrm{class}' \restriction U'$, where $f, f' \in \mathscr{A}$. Then, $\psi = f'^{-1} \circ f$, defined on $f^{-1}(f'(U'))$ to $U'$, is a local diffeomorphism that maps $\xi = f(\mathrm{class}(r))$ to $\xi' = f'(\mathrm{class}'(r'))$.

Then, according to (§4), $n = n'$ and there exists a local diffeomorphism $\phi$ of $\mathbf{R}^n$, lifting locally $\psi$ and mapping $r$ to $r'$. Its germ realizes an arrow of the groupoid **G** connecting (F, $r$) to (F', $r'$). Of course, when F($r$) $\neq$ F'($r'$) there cannot be an arrow, by definition. Therefore, as in the more restrictive case of orbifolds, the fibers of the evaluation map are the transitive components of the structure groupoid **G** of the quasifold.

Now, the theorem follows the formal flow of (*op. cit.* §5): let $\mathscr{A}$ and $\mathscr{A}'$ be two atlases of X and consider $\mathscr{A}'' = \mathscr{A} \coprod \mathscr{A}'$. With an obvious choice of notation: Obj(**G**″) = Obj(**G**) $\coprod$ Obj(**G**′) and **G**″ contains naturally **G** and **G**′ as full subgroupoids. The question then is: how does the adjunction of the crossed arrows between **G** and **G**′ change the distribution of transitivity-components? According to the previous proposition, it changes nothing since, for **G**, **G**′ or **G**″, the set of transitivity-components are always exactly the fibers of the respective subductions ev. In other words, the set of groupoid components is always X, for any atlas of X. Thus **G** and **G**′ are equivalent to **G**″, therefore **G** and **G**′ are equivalent. $\square$

7. GENERAL DESCRIPTION OF THE STRUCTURE GROUPOID. — The general description of the structure groupoid of a quasifold X follows exactly the description in the case of orbifolds (*op. cit.*). We remind it here for clarity. Let X be a quasifold. Let $\mathscr{A}$ be an atlas, let $\mathscr{F}$ be the associated strict generating family, and let **G** be the associated groupoid. We know from the previous paragraph that the groupoid components in Obj(**G**) are the fibers of the projection ev: (F, $r$) $\mapsto$ F($r$). Then, the (algebraic) structure of the groupoid reduces to the algebraic structure of each full subgroupoid **G**$_x$, $x \in$ X, that is,

$$\begin{cases} \text{Obj}(\mathbf{G}_x) &= \{(F, r) \in \mathcal{N} \mid F(r) = x\}, \\ \text{Mor}(\mathbf{G}_x) &= \{\mathbf{g} \in \text{Mor}_{\mathbf{G}}((F, r), (F', r')) \mid F(r) = x\}; \end{cases}$$

more precisely, $\mathbf{g} = \text{germ}(\phi)_r$ where $\phi$ is a local diffeomorphism defined in the domain of F to the domain of F′, mapping $r$ to $r'$ and such that F′ $\circ$ $\phi$ =$_{\text{loc}}$ F on an open neighborhood of $r$. In other words,

$$\text{Obj}(\mathbf{G}_x) = \text{ev}^{-1}(x) \quad \text{and} \quad \text{Mor}(\mathbf{G}_x) = (\text{ev} \circ \text{src})^{-1}(x).$$

Let $f$ be a chart in $\mathscr{A}$, let U = dom($f$) and let $\tilde{\text{U}}$ = class$^{-1}$(U) $\subset \mathbf{R}^n$ be the domain of its strict lifting F = $f \circ$ class $\upharpoonright \tilde{\text{U}}$, where class: $\mathbf{R}^n \to \mathbf{R}^n/\Gamma$. Without loss of generality, we shall assume that the domains of all charts, and thus the domains of the strict liftings, are connected.

The subgroupoid **G**$_x$ is the assemblage of the subgroupoids **G**$_x^{\text{F}}$. For all F $\in \mathscr{F}$,

$$\begin{cases} \text{Obj}(\mathbf{G}_x^{\text{F}}) &= \{\text{F}\} \times \text{dom}(\text{F}), \\ \text{Mor}(\mathbf{G}_x^{\text{F}}) &= \{\text{germ}(\phi)_r \in \text{Mor}(\mathbf{G}_x) \mid r, \phi(r) \in \text{dom}(\text{F})\}. \end{cases}$$

That is, Mor(**G**$_x^{\text{F}}$) = src$^{-1}$(Obj(**G**$_x^{\text{F}}$)) $\cap$ trg$^{-1}$(Obj(**G**$_x^{\text{F}}$)). The assemblage is made first by connecting the groupoid **G**$_x^{\text{F}}$ to **G**$_x^{\text{F}'}$ with any arrow germ($\phi$)$_r$, from (F, $r$) to (F′, $r'$) such that $x = $ F($r$) = F($r'$) and $\phi(r) = r'$. Secondly, by spreading the arrows by composition. We can represent this construction by a groupoid-set-theoretical diagram:

$$\mathbf{G} = \coprod_{x \in X} \mathbf{G}_x \quad \text{and} \quad \mathbf{G}_x = \mathbf{G}_x^{F_1} - \mathbf{G}_x^{F_2} - \cdots - \mathbf{G}_x^{F_{N_x}}$$

where the $F_i$'s are the charts having $x$ in their images and $N_x$ is the number of such charts (the atlas $\mathscr{A}$ is assumed locally finite). The link between two groupoids: $\mathbf{G}_x^{F_i} - \mathbf{G}_x^{F_j}$ represents the spreading of the arrows by adjunction of one of them. Note that this is absolutely not a smooth representation of $\mathbf{G}$, since the projection $\mathrm{ev} \circ \mathrm{src} \colon \mathrm{Mor}(\mathbf{G}) \to X$ is a subduction. Moreover, the order of assembly has no influence on the result.

EXAMPLES. In the case of orbifolds, where the structure group is finite, this assemblage of groupoids can be completely visual: for example, the teardrop in [IZL17, Figure 3]. It is more difficult in the case of a strict quasifold, with dense structure group. For example, the irrational torus $T_\alpha = \mathbf{R}/(\mathbf{Z} + \alpha \mathbf{Z})$, which was described as a diffeological space in [PDPI83] for the first time. Now, with the identification of this new subcategory {Quasifolds} in {Diffeology}, the irrational torus becomes a *quasitorus*.[9] For the generating family $\{\mathrm{class}\colon \mathbf{R} \to T_\alpha\}$, the objects of the structure groupoid equal just $\mathbf{R}$. Moreover, in this simple case, as we see in (§11), the groupoid $\mathbf{G}_\alpha$ is the groupoid of the action of the subgroup $\mathbf{Z} + \alpha \mathbf{Z}$ by translation. Therefore, one has

$$\mathrm{Obj}(\mathbf{G}_\alpha) = \mathbf{R} \quad \text{and} \quad \mathrm{Mor}(\mathbf{G}_\alpha) = \{(x, \mathbf{t}_{n+\alpha m}) \mid x \in \mathbf{R} \text{ and } n, m \in \mathbf{Z}\},$$

where the bold letter $\mathbf{t}$ denotes a translation. The source and target are given by

$$\mathrm{src}(x, \mathbf{t}_{n+\alpha m}) = x \quad \text{and} \quad \mathrm{trg}(x, \mathbf{t}_{n+\alpha m}) = x + n + \alpha m.$$

Also, the composition of arrows is given by

$$(x, \mathbf{t}_{n+\alpha m}) \cdot (x + n + \alpha m, \mathbf{t}_{n'+\alpha m'}) = (x, \mathbf{t}_{n+n'+\alpha(m+m')}).$$

The subgroupoid $\mathbf{G}_{\alpha,\tau}^{\mathrm{class}}$, with $\tau \in T_\alpha$, is then

$$\mathbf{G}_{\alpha,\tau}^{\mathrm{class}} = \{(x, \mathbf{t}_k) \mid \mathrm{class}(x) = \tau \text{ and } k \in \mathbf{Z} + \alpha \mathbf{Z}\}.$$

For example, at $\tau = 0$ we get $\mathbf{G}_{\alpha,0}^{\mathrm{class}} = \{(n + \alpha m, \mathbf{t}_{n'+\alpha m'}) \mid n, n', m, m' \in \mathbf{Z}\}$.

8. THE STRUCTURE GROUPOID IS ÉTALE AND HAUSDORFF. — Let $\mathscr{A}$ be an atlas of a quasifold X. The structure groupoid $\mathbf{G}$ associated with the generating family of the atlas $\mathscr{A}$ is étale, namely: the projection $\mathrm{src}\colon \mathrm{Mor}(\mathbf{G}) \to \mathrm{Obj}(\mathbf{G})$ is an étale smooth map, that is, a local diffeomorphism at each point [PIZ13, §2.5].

PROPOSITION 1. *For all $\mathbf{g} \in \mathrm{Mor}(\mathbf{G})$, there exists a D-open superset $\mathscr{O}$ of $\mathbf{g}$ such that* $\mathrm{src}$ *restricted to $\mathscr{O}$ is a local diffeomorphism.*

PROPOSITION 2. *The groupoid $\mathbf{G}$ is locally compact and Hausdorff.*

NOTE. Since the atlas $\mathscr{A}$ is assumed to be locally finite, the preimages of the objects of $\mathbf{G}$ by the source map, or the target map, are countable.

---

[9]Or, in this case, a quasicircle.

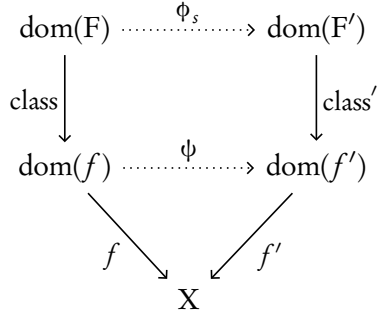

FIGURE 3. Lifting local diffeomorphisms

*Proof.* This proof is the same as in the case of orbifolds [IZL17, §7]. We just have to pay attention to the fact the structure group is now countable, and not just finite.

(1) Let us first check that the groupoid **G** is étale. That is, src: $\mathrm{Mor}(\mathbf{G}) \to \mathrm{Obj}(\mathbf{G})$ is everywhere a local diffeomorphism.

Let us pick a germ $\mathbf{g} = \mathrm{germ}(\Phi)_\nu \in \mathrm{Mor}(\mathbf{G})$, with $\nu = \mathrm{src}(\mathbf{g}) = (\mathrm{F}, r)$ and $\mathrm{trg}(\mathbf{g}) = (\mathrm{F}', r')$. Thus, $\Phi$ is defined by some $\phi \in \mathrm{Diff}_{\mathrm{loc}}(\mathbf{R}^n)$ with $\mathrm{dom}(\phi) \subset \mathrm{dom}(\mathrm{F})$, $r' = \phi(r) \in \mathrm{dom}(\mathrm{F}')$ and such that $\mathrm{F}' \circ \phi = \mathrm{F} \restriction \mathscr{B}$. We choose $\phi \colon \mathscr{B} \to \mathrm{dom}(\mathrm{F}')$ to be defined on a small ball centered at $r$. By abuse of notation we shall denote $\mathbf{g} = \mathrm{germ}(\phi)_r$, where $\phi \in \mathrm{Diff}_{\mathrm{loc}}(\mathrm{dom}(\mathrm{F}), \mathrm{dom}(\mathrm{F}'))$. That is, $\phi$ now contains implicitly the data source and target. Now, let

$$\mathrm{F} = f \circ \mathrm{class} \quad \text{and} \quad \mathrm{F}' = f' \circ \mathrm{class}',$$

where $f$ and $f'$ belong to $\mathscr{A}$, class: $\mathbf{R}^n \to \mathbf{R}^n/\Gamma$ and class$'$: $\mathbf{R}^n \to \mathbf{R}^n/\Gamma'$ are the projections. If $\psi$ is the transition map $f'^{-1} \circ f$, then $\mathrm{class}(r) \in \mathrm{dom}(\psi)$ and $\psi(\mathrm{class}(r)) = \mathrm{class}'(r')$. This situation is illustrated by the diagram of Figure 3, where, except for the family $\phi_s$ which will vary around $\phi$, the vertices and arrows are fixed as soon as the representant $\Phi$ of the germ $\mathbf{g}$ is chosen. Now, let

$$\mathscr{O} = \{\mathrm{germ}(\phi)_t \mid t \in \mathscr{B}\} \subset \mathrm{Mor}(\mathbf{G}).$$

Hence, src $\restriction \mathscr{O}\colon \mathrm{germ}(\phi)_t \mapsto t$ is smooth and injective,[10] as well as its inverse $t \mapsto \mathrm{germ}(\phi)_t$, which is defined on $\mathscr{B}$. Let us now show that $\mathscr{O}$ is a D-open subset, That is, for each plot $\mathrm{P}\colon s \mapsto \mathbf{g}_s$ in $\mathrm{Mor}(\mathbf{G})$, the subset $\mathrm{P}^{-1}(\mathscr{O}) \subset \mathrm{dom}(\mathrm{P})$ is open. Let $s \in \mathrm{P}^{-1}(\mathscr{O})$, that is, $\mathbf{g}_s \in \mathscr{O}$, *i.e.* $\mathbf{g}_s = \mathrm{germ}(\phi)_{r_s}$, where $r_s = \mathrm{src}(\mathbf{g}_s)$, the discrete index F here is implicit.

Then, for all $s \in \mathrm{dom}(\mathrm{P})$, there exists a small ball $\mathscr{V}$ centered at $s$ and a plot $s' \mapsto (\phi_{s'}, r_{s'})$, defined on $\mathscr{V}$, such that $\mathbf{g}_{s'} = \mathrm{germ}(\phi_{s'})_{r_{s'}}$ with $\mathrm{germ}(\phi_s)_{r_s} = \mathrm{germ}(\phi)_{r_s}$ and $r_{s'} \in \mathscr{B}$.

---

[10]Maybe we should recall that $\mathrm{germ}(\phi)_t = \mathrm{germ}(\phi')_{t'}$ if and only if: $t = t'$ and there exists an open ball $\mathscr{B}$ centered at $t$ such that $\phi \restriction \mathscr{B} = \phi' \restriction \mathscr{B}$.

Since $s' \mapsto \phi_{s'}$ is smooth, by definition the subset

$$\{(s', r) \in \mathcal{V} \times \mathcal{B} \mid r \in \mathrm{dom}(\phi_{s'})\}$$

is necessarily open. Since it contains $(s, r_s)$, it contains a product $\mathcal{V}' \times \mathcal{B}'$, where $\mathcal{V}'$ is a small ball centered at $s$ and $\mathcal{B}'$ is a small ball centered at $r_s$. This implies that, for all $s' \in \mathcal{V}'$, $\mathcal{B}' \subset \mathrm{dom}(\phi_{s'})$. In particular, $\mathcal{B}' \subset \mathrm{dom}(\phi)$.

Then, $\xi_s = \phi_s \circ \phi^{-1} \colon \phi(\mathcal{B}') \to \mathrm{dom}(F')$ is a local diffeomorphism of $\mathrm{dom}(F')$. However, for all $s'$, one has $\mathrm{class}' \circ \phi_{s'} = \psi \circ \mathrm{class}$, wherever it is defined. This is shown by the above diagram, where the dots denote a local map.

Thus, $\mathrm{class}' \circ \xi_{s'} = \mathrm{class}'$. Indeed, $\mathrm{class}' \circ \xi_{s'} = \mathrm{class}' \circ \phi_{s'} \circ \phi^{-1} = \psi \circ \mathrm{class} \circ \phi^{-1}$ and $\psi \circ \mathrm{class} = \mathrm{class}' \circ \phi$. Now, thanks to (§7), for all $s' \in \mathcal{V}'$ there is a $\gamma' \in \Gamma'$ such that $\xi_{s'} = \gamma'$, and the map $s' \mapsto \gamma'$ is smooth. Then, since $\mathcal{V}'$ is connected and $\Gamma'$ is discrete [PIZ13, Exercise 8], $\gamma'$ is constant on $\mathcal{V}'$. Now $s \in \mathcal{V}'$, thus, for $s' = s$, $\gamma' = \phi^{-1} \circ \phi_s = \phi^{-1} \circ \phi = \mathbf{1}$. Hence, $\phi_{s'} = \phi$ on $\mathcal{V}'$, and $\mathbf{g}_{s'} = \mathrm{germ}(\phi)_{r_{s'}}$ on $\mathcal{V}'$, that is, $P(\mathcal{V}') \subset \mathcal{O}$. Then, each $s \in \mathrm{dom}(P)$ such that $P(s) \in \mathcal{O}$ is the center of an open ball whose image is contained in $\mathcal{O}$. Therefore, $P^{-1}(\mathcal{O})$ is a union of open balls, thus $P^{-1}(\mathcal{O})$ is open and $\mathcal{O}$ is D-open. Thus, the map $\mathrm{src} \colon \mathrm{germ}(\phi)_t \mapsto t$, restricted to $\mathcal{O}$, is a local diffeomorphism: the source map is étale.

(2) Next, let us check that $\mathrm{Mor}(\mathbf{G})$ is Hausdorff. As above, let $\mathbf{g} = \mathrm{germ}(\phi)_r \in \mathrm{Mor}(\mathbf{G})$. We can also represent $\mathbf{g}$ by a triple $(F, r, \mathrm{germ}(\phi)_r)$, with $\phi \in \mathrm{Diff}_{\mathrm{loc}}(\mathrm{dom}(F), \mathrm{dom}(F'))$. Then, let $\mathbf{g}' = \mathrm{germ}(\psi)_s$ be another germ represented by $(G, s, \mathrm{germ}(\psi)_s)$, different from $\mathbf{g}$, with $\psi \in \mathrm{Diff}_{\mathrm{loc}}(\mathrm{dom}(G), \mathrm{dom}(G'))$. We separate the situation in three cases:

$$F \neq G$$
$$\text{or}$$
$$F = G \quad \begin{cases} r \neq s \\ \text{or} \\ r = s \quad \text{but} \quad \mathrm{germ}(\phi)_r \neq \mathrm{germ}(\psi)_s \end{cases}$$

In the first two cases ($F \neq G$, and $F = G$ but $r \neq s$), since the source map is étale and since the Nebula is Hausdorff, it is sufficient to consider two small separated balls $\mathcal{B}$ and $\mathcal{B}'$, centered around $r$ and $s$, to get two D-open subsets of $\mathrm{Mor}(\mathbf{G})$ that separate the two different germs. Indeed, let $\mathcal{O} = \mathrm{src}^{-1}(\mathcal{B})$ and $\mathcal{O}' = \mathrm{src}^{-1}(\mathcal{B}')$ be the D-open subset on which the source map is a local diffeomorphism. If there were a point $\mathbf{g}'' \in \mathcal{O} \cap \mathcal{O}'$, then $\mathrm{src}(\mathbf{g}'')$ would belong to $\mathcal{B} \cap \mathcal{B}'$, which is empty.

The last case ($(F, r, \mathrm{germ}(\phi)_r)$ and $(F, r, \mathrm{germ}(\psi)_r)$ with $\mathrm{germ}(\phi)_r \neq \mathrm{germ}(\psi)_r$) divides in two subcases: $\mathrm{codom}(\phi) \neq \mathrm{codom}(\psi)$ and $\mathrm{codom}(\phi) = \mathrm{codom}(\psi)$.

In the first sub-case, when $\mathrm{codom}(\phi) \neq \mathrm{codom}(\psi)$, since the codomains are different, we consider a small ball $\mathcal{B}$ around $r$ such that its images by $\phi$ and $\psi$ are separated. Then $\mathcal{O} = \{\mathrm{germ}(\phi)_t \mid t \in \mathcal{B}\}$ and $\mathcal{O}' = \{\mathrm{germ}(\psi)_t \mid t \in \mathcal{B}\}$ are two open subsets in $\mathrm{Mor}(\mathbf{G})$ that separate $\mathbf{g}$ and $\mathbf{g}'$, since no germ in $\mathcal{O}'$ has the same codomain as any germ in $\mathcal{O}$.

In the second sub-case, when $\mathrm{codom}(\phi) = \mathrm{codom}(\psi)$, let us consider the composite $f = \phi \circ \psi^{-1}$, defined on an open neighborhood of $\psi(r)$. Thanks to the theorem of (§3), $f(s) =_{\mathrm{loc}} \gamma' \cdot s$, for some $\gamma' \in \Gamma'$, which is the structure group of the quasifold for the plot $F'$. Since we have assumed that $\mathrm{germ}(\phi)_r \neq \mathrm{germ}(\psi)_r$, we have that $\gamma' \neq \mathbf{1}$. Hence, there is a small ball $\mathscr{B}$ around $r$ on which $\phi = \gamma' \circ \psi$. Let $\mathscr{O}' = \{(F, t, \mathrm{germ}(\psi)_t) \mid t \in \mathscr{B}\}$ and $\mathscr{O} = \{(F, t, \mathrm{germ}(\phi)_t) \mid t \in \mathscr{B}\}$, that is, $\mathscr{O} = \{(F, t, \gamma' \circ \mathrm{germ}(\psi)_t) \mid t \in \mathscr{B}\}$. As we know, they are two D-open subsets in $\mathrm{Mor}(\mathbf{G})$ and, since $\gamma' \neq \mathbf{1}$, we have that $\mathscr{O} \cap \mathscr{O}' = \emptyset$. Therefore, the two germs $(F, r, \mathrm{germ}(\phi)_r)$ and $(F, r, \mathrm{germ}(\psi)_r)$ are separated.

In conclusion, $\mathrm{Mor}(\mathbf{G})$ is Hausdorff for the D-topology.

According to the previous Note, the preimages of an object $(F, r) \in \mathcal{N} = \mathrm{Obj}(\mathbf{G})$ are the germs of all the local diffeomorphisms $\Phi \colon (F, r) \mapsto (F', r')$, such that $F(r) = r'$ and $F =_{\mathrm{loc}} F' \circ \phi$ around $r$, where $\phi$ is a local diffeomorphism of $\mathbf{R}^n$. Since the atlas $\mathscr{A}$ is locally finite, there are a finite number of $F' \in \mathscr{F}$ such that $F(r) = F'(r')$. Now, for such $F'$ the number of $r' \in \mathrm{dom}\, F'$ such that $F'(r') = F(r)$ is at most equal to the number of elements of the structure group $\Gamma'$, that is countable. Therefore, the preimages of $(F, r)$ by the source map is countable, and that works obviously in the same way for the preimages of the target map. □

9. MRW-EQUIVALENCE OF STRUCTURE GROUPOIDS. — We consider a quasifold X and two atlases $\mathscr{A}$ and $\mathscr{A}'$, with associated strict generating families $\mathscr{F}$ and $\mathscr{F}'$. We shall show in this section that the associated groupoids are equivalent in the sense of Muhly-Renault-Williams [MRW87, 2.1]; this will later give Morita-equivalent $\mathbf{C}^*$-algebras.

This section follows [IZL17, §8]; we just check that the fact that the structure groups are countable and not just finite, does not change the result.

Let us recall what is an MRW-equivalence of groupoids. Let $\mathbf{G}$ and $\mathbf{G}'$ be two locally compact groupoids. We say that a locally compact space Z is a $(\mathbf{G}, \mathbf{G}')$-equivalence if

(i) Z is a left principal $\mathbf{G}$-space.
(ii) Z is a right principal $\mathbf{G}'$-space.
(iii) The $\mathbf{G}$ and $\mathbf{G}'$ actions commute.
(iv) The action of $\mathbf{G}$ on Z induces a bijection of $\mathrm{Z}/\mathbf{G}$ onto $\mathrm{Obj}(\mathbf{G}')$.
(v) The action of $\mathbf{G}'$ on Z induces a bijection of $\mathrm{Z}/\mathbf{G}'$ onto $\mathrm{Obj}(\mathbf{G})$.

Let $\mathrm{src}\colon \mathrm{Z} \to \mathrm{Obj}(\mathbf{G})$ and $\mathrm{trg}\colon \mathrm{Z} \to \mathrm{Obj}(\mathbf{G}')$ be the maps defining the composable pairs associated with the actions of $\mathbf{G}$ and $\mathbf{G}'$. That is, a pair $(\mathbf{g}, z)$ is composable if $\mathrm{trg}(\mathbf{g}) = \mathrm{src}(z)$, and the composite is denoted by $\mathbf{g} \cdot z$. Moreover, a pair $(\mathbf{g}', z)$ is composable if $\mathrm{src}(\mathbf{g}') = \mathrm{trg}(z)$, and the composite is denoted by $z \cdot \mathbf{g}'$.

Let us also recall that an action is principal in the sense of Muhly-Renault-Williams, if it is free: $\mathbf{g} \cdot z = z$ only if $\mathbf{g}$ is a unit, and the *action map* $(\mathbf{g}, z) \mapsto (\mathbf{g} \cdot z, z)$, defined on the composable pairs, is proper [MRW87, §2].

Now, using the hypothesis and notations of (§6), let us define Z to be the space of germs of local diffeomorphisms, from the nebula of the family $\mathscr{F}$ to the nebula of the family

$\mathscr{F}'$, that project on the identity by the evaluation map. That is,

$$Z = \left\{ \operatorname{germ}(f)_r \;\middle|\; \begin{array}{l} f \in \operatorname{Diff}_{\operatorname{loc}}(\operatorname{dom}(F), \operatorname{dom}(F')), r \in \operatorname{dom}(F), \\ F \in \mathscr{F}, F' \in \mathscr{F}' \text{ and } F' \circ f = F \upharpoonright \operatorname{dom}(f). \end{array} \right\}.$$

Let[11]

$$\operatorname{src}(\operatorname{germ}(f)_r) = r \quad \text{and} \quad \operatorname{trg}(\operatorname{germ}(f)_r) = f(r).$$

Then, the action of $\mathbf{g} \in \operatorname{Mor}(\mathbf{G})$ on $\operatorname{germ}(f)_r$ is defined by composition if $\operatorname{trg}(\mathbf{g}) = r$, that is, $\mathbf{g} \cdot \operatorname{germ}(f)_r = \operatorname{germ}(f \circ \phi)_s$, where $\mathbf{g} = \operatorname{germ}(\phi)_s$, $\phi \in \operatorname{Diff}_{\operatorname{loc}}(\mathcal{N})$ and $\phi(s) = r$. Symmetrically, the action of $\mathbf{g}' \in \operatorname{Mor}(\mathbf{G}')$ on $\operatorname{germ}(f)_r$ is defined if $\operatorname{src}(\mathbf{g}') = f(r)$ by $z \cdot \mathbf{g}' = \operatorname{germ}(\phi' \circ f)_r$, where $\mathbf{g}' = \operatorname{germ}(\phi')_{f(r)}$. Then, we have:

THEOREM. *The actions of $\mathbf{G}$ and $\mathbf{G}'$ on $Z$ are principal, and $Z$ is a $(\mathbf{G}, \mathbf{G}')$-equivalence in the sense of Muhly-Renault-Williams.*

*Proof.* First of all, let us point out that $Z$ is a subspace of the morphisms of the groupoid $\mathbf{G}''$ built in (§6) by adjunction of $\mathbf{G}$ and $\mathbf{G}'$, and is equipped with the subset diffeology. All these groupoids are locally compact and Hausdorff (§8).

Let us check that the action of $\mathbf{G}$ on $Z$ is free. In our case, $z = \operatorname{germ}(f)_r$ and $\mathbf{g} = \operatorname{germ}(\phi)_s$, where $f$ and $\phi$ are local diffeomorphisms. If $\mathbf{g} \cdot z = z$, then obviously $\mathbf{g} = \operatorname{germ}(\mathbf{1})_r$.

Next, let us denote by $\rho$ the action of $\mathbf{G}$ on $Z$, defined on

$$\mathbf{G} \star Z = \{(\mathbf{g}, z) \in \operatorname{Mor}(\mathbf{G}) \times Z \mid \operatorname{trg}(\mathbf{g}) = \operatorname{src}(z)\} \quad \text{by} \quad \rho(\mathbf{g}, z) = \mathbf{g} \cdot z.$$

This action is smooth because the composition of local diffeomorphisms is smooth, and passes onto the quotient groupoid in a smooth operation, see [IZL17, §3]. Moreover, this action is invertible, its inverse being defined on

$$Z \star Z = \{(z', z) \in Z \times Z \mid \operatorname{trg}(z') = \operatorname{trg}(z)\} \quad \text{by} \quad \rho^{-1}(z', z) = (\mathbf{g} = z' \cdot z^{-1}, z).$$

In detail, $\rho^{-1}(\operatorname{germ}(h)_s, \operatorname{germ}(f)_r) = (\operatorname{germ}(f^{-1} \circ h)_s, \operatorname{germ}(f)_r)$, with $f(r) = h(s)$. Now, the inverse is also smooth, when $Z \star Z \subset Z \times Z$ is equipped with the subset diffeology. In other words, $\rho$ is an induction, that is, a diffeomorphism from $\mathbf{G} \star Z$ to $Z \star Z$. However, since $\mathbf{G} \star Z$ and $Z \star Z$ are defined by closed relations, and $\mathbf{G}$ and $Z$ are Hausdorff, $\mathbf{G} \star Z$ and $Z \star Z$ are closed into their ambient spaces. Thus, the intersection of a compact subset in $Z \times Z$ with $Z \star Z$ is compact, and its preimage by the induction $\rho$ is compact. Therefore, $\rho$ is proper. We notice that the fact that the structure groups are no longer finite but just countable does not play a role here.

It remains to check that the action of $\mathbf{G}$ on $Z$ induces a bijection of $Z/\mathbf{G}$ onto $\operatorname{Obj}(\mathbf{G}')$. Let us consider the map class: $Z \to \operatorname{Obj}(\mathbf{G}')$ defined by $\operatorname{class}(\operatorname{germ}(f)_r) = f(r)$. Then, let $\operatorname{class}(z) = \operatorname{class}(z')$, with $z = \operatorname{germ}(f)_r$ and $z' = \operatorname{germ}(f')_{r'}$, that is, $f(r) = f'(r')$. However, since $f$ and $f'$ are local diffeomorphisms, $\phi = f'^{-1} \circ f$ is a local diffeomorphism with $\phi(r') = r$. Let $\mathbf{g} = \mathrm{germ}(\phi)_{r'}$, then $\mathbf{g} \in \mathrm{Mor}(\mathbf{G})$ and $z' = \mathbf{g} \cdot z$. Hence, the map class projects onto an injection from $Z/\mathbf{G}$ to $\mathrm{Obj}(\mathbf{G}')$. Now, let $(F', r') \in \mathrm{Obj}(\mathbf{G}')$, and let $x = F'(r') \in X$. Since $\mathscr{F}$ is a generating family, there exist $(F, r) \in \mathrm{Obj}(\mathbf{G})$ such that $F(r) = x$. Let $\psi$ and $\psi'$ be the charts of X defined by factorization: $F = \psi \circ \mathrm{class}$ and $F' = \psi' \circ \mathrm{class}'$, where $\mathrm{class}: \mathbf{R}^n \to \mathbf{R}^n/\Gamma$ and $\mathrm{class}': \mathbf{R}^n \to \mathbf{R}^n/\Gamma'$. Let $\xi = \mathrm{class}(r)$ and $\xi' = \mathrm{class}'(r')$. Since $\psi(\xi) = \psi'(\xi') = x$, $\Psi =_{\mathrm{loc}} \psi'^{-1} \circ \psi$ is a local diffeomorphism from $\mathbf{R}^n/\Gamma$ to $\mathbf{R}^n/\Gamma'$ mapping $\xi$ to $\xi'$. Hence, according to (§4), there exists a local diffeomorphism $f$ from $\mathrm{dom}(F)$ to $\mathrm{dom}(F')$, such that $\mathrm{class}' \circ f = \Psi \circ \mathrm{class}$ and $f(r) = r'$. Thus, $z = \mathrm{germ}(f)_r$ belongs to Z and $\mathrm{class}(z) = r'$ (precisely the element $(F', r')$ of the nebula of $\mathscr{F}'$). Therefore, the injective map class from $Z/\mathbf{G}$ to $\mathrm{Obj}(\mathbf{G}')$ is also surjective, and identifies the two spaces. Obviously, what has been said for the side $\mathbf{G}$ can be translated to the side $\mathbf{G}'$; the construction is completely symmetric. In conclusion, Z satisfies the conditions of a $(\mathbf{G}, \mathbf{G}')$-equivalence, in the sense of Muhly-Renault-Williams. □

---

[11] For the sake of simplicity, we make an abuse of notation: in reality one should write, more precisely, $\operatorname{src}(\operatorname{germ}(f)_r) = (F, r)$ and $\operatorname{trg}(\operatorname{germ}(f)_r) = (F', f(r))$.

## The $\mathbf{C}^*$-Algebras Of A Quasifold

We use the construction of the $\mathbf{C}^*$-Algebra associated with an arbitrary locally compact groupoid $\mathbf{G}$, equipped with a Haar system, introduced and described by Jean Renault in [JR80, Part II, §1]. Note that, for this construction, only the topology of the groupoid is involved, and diffeological groupoids, when regarded as topological groupoids, are equipped with the D-topology[12].

We will denote by $\mathscr{C}(\mathbf{G})$ the completion of the compactly supported continuous complex functions on $\mathrm{Mor}(\mathbf{G})$, for the uniform norm. And we still consider, as is done for orbifolds, the particular case where the Haar system is given by the *counting measure*. Let $f$ and $g$ be two compactly supported complex functions, the convolution and the involution are defined by

$$f * g(\gamma) = \sum_{\beta \in \mathbf{G}^x} f(\beta \cdot \gamma) g(\beta^{-1}) \quad \text{and} \quad f^*(\gamma) = f(\gamma^{-1})^*.$$

The sums involved are supposed to converge. Here, $\gamma \in \mathrm{Mor}(\mathbf{G})$, $x = \mathrm{src}(\gamma)$ and $\mathbf{G}^x = \mathrm{trg}^{-1}(x)$ is the subset of arrows with target $x$. The star in $z^*$ denotes the conjugate of the complex number $z$. By definition, the vector space $\mathscr{C}(\mathbf{G})$, equipped with these two operations, is the $\mathbf{C}^*$-algebra associated with the groupoid $\mathbf{G}$.

10. The $\mathbf{C}^*$-algebra of a quasifold. — Let X be a quasifold, let $\mathscr{A}$ be an atlas and let $\mathbf{G}$ be the structure groupoid associated with $\mathscr{A}$. Since the atlas $\mathscr{A}$ is locally finite, the convolution defined above is well defined. Indeed, in this case:

Proposition. *For every compactly supported complex function $f$ on $\mathbf{G}$, for all $\nu = (F, r) \in \mathscr{N} = \mathrm{Obj}(\mathbf{G})$, the set of arrows $\mathbf{g} \in \mathbf{G}^\nu$ such that $f(\mathbf{g}) \neq 0$ is finite. That is, $\#\mathrm{Supp}(f \upharpoonright \mathbf{G}^\nu) < \infty$. The convolution is then well defined on $\mathscr{C}(\mathbf{G})$.*

---

[12]Since smooth maps are D-continuous and diffeomorphism are D-homeomorphisms.

Then, for each atlas $\mathscr{A}$ of the quasifold X, we get the $\mathbf{C}^*$-algebra $\mathfrak{A} = (\mathscr{C}(\mathbf{G}), *)$. The dependence of the $\mathbf{C}^*$-algebra on the atlas is given by the following theorem, which is a generalization of [IZL17, §9].

THEOREM. *Different atlases give Morita-equivalent $\mathbf{C}^*$-algebras. Diffeomorphic quasifolds have Morita-equivalent $\mathbf{C}^*$-algebras.*

In other words, we have defined a functor from the subcategory of isomorphic {Quasifolds} in diffeology, to the category of Morita-equivalent {$\mathbf{C}^*$-Algebras}.

*Proof.* Considering the proposition, $\mathbf{G}^\nu = \mathrm{trg}^{-1}(\nu)$ with $\nu \in \mathrm{Obj}(\mathbf{G})$. The space of objects of $\mathbf{G}$ is a disjoint sum of Euclidean domains, thus $\{\nu\}$ is a closed subset. Now, $\mathrm{trg}\colon \mathrm{Mor}(\mathbf{G}) \to \mathrm{Obj}(\mathbf{G})$ is smooth then continuous, for the D-topology. Hence, $\mathbf{G}^\nu = \mathrm{trg}^{-1}(\nu)$ is closed and countable by (§8). Now, $\mathrm{Supp}(f \restriction \mathbf{G}^\nu) = \mathrm{Supp}(f) \cap \mathbf{G}^\nu$ is the intersection of a compact and a closed countable subspace, thus it is compact and countable, that is finite.

Next, thanks to (§9), different atlases give equivalent groupoids in the sense of Muhly-Renault-Williams. Moreover, thanks to [MRW87, Thm. 2.8], different atlases give strongly Morita-equivalent $\mathbf{C}^*$-algebras. Therefore, diffeomorphic quasifolds have strongly Morita-equivalent associated $\mathbf{C}^*$-algebras. □

11. THE $\mathbf{C}^*$-ALGEBRA OF THE IRRATIONAL TORUS. — The first and most famous example is the so-called Denjoy-Poincaré torus, or irrational torus, or noncommutative torus, or, more recently, quasitorus. It is, according to its first definition, the quotient set of the 2-torus $\mathrm{T}^2$ by the irrational flow of slope $\alpha \in \mathbf{R} - \mathbf{Q}$. We denote it by $\mathrm{T}_\alpha = \mathrm{T}^2/\Delta_\alpha$, where $\Delta_\alpha$ is the image of the line $y = \alpha x$ by the projection $\mathbf{R}^2 \to \mathrm{T}^2 = \mathbf{R}^2/\mathbf{Z}^2$. This space has been the first example studied with the tools of diffeology, in [PDPI83], where many non trivial properties have been highlighted.[13] Diffeologically speaking,

$$\mathrm{T}_\alpha \simeq \mathbf{R}/(\mathbf{Z} + \alpha \mathbf{Z}).$$

The composite

$$\mathbf{R} \xrightarrow{\mathrm{class}} \mathbf{R}/(\mathbf{Z} + \alpha \mathbf{Z}) \xrightarrow{f} \mathrm{T}_\alpha, \text{ with } \mathrm{F} = f \circ \mathrm{class},$$

summarizes the situation where $\mathscr{A} = \{f \colon \mathbf{R}/(\mathbf{Z} + \alpha \mathbf{Z}) \to \mathrm{T}_\alpha\}$ is the canonical atlas of $\mathrm{T}_\alpha$, containing the only chart $f$, and $\mathscr{F} = \{\mathrm{F} = f \circ \mathrm{class}\}$ is the associated canonical strict generating family. According to the above (§3), the groupoid $\mathbf{G}_\alpha$ associated with the atlas $\mathscr{A}$ is simply

$$\mathrm{Obj}(\mathbf{G}_\alpha) = \mathbf{R} \quad \text{and} \quad \mathrm{Mor}(\mathbf{G}_\alpha) = \{(x, \mathbf{t}_{n+\alpha m}) \mid x \in \mathbf{R} \text{ and } n, m \in \mathbf{Z}\}.$$

However, we can also identify $\mathrm{T}_\alpha$ with $(\mathbf{R}/\mathbf{Z})/[(\mathbf{Z} + \alpha \mathbf{Z})/\mathbf{Z}]$, that is

$$\mathrm{T}_\alpha \simeq \mathrm{S}^1/\mathbf{Z}, \text{ with } \quad \underline{m}(z) = e^{2i\pi\alpha m} z,$$

---

[13]See for example Exercise 4 and §8.39 in [PIZ13].

for all $m \in \mathbf{Z}$ and $z \in S^1$. Moreover, the groupoid $\mathbf{S}$ of this action of $\mathbf{Z}$ on $S^1 \subset \mathbf{C}$ is simply

$$\mathrm{Obj}(\mathbf{S}_\alpha) = S^1 \quad \text{and} \quad \mathrm{Mor}(\mathbf{S}_\alpha) = \{(z, e^{2i\pi\alpha m}) \mid z \in S^1 \text{ and } m \in \mathbf{Z}\}.$$

The groupoids $\mathbf{G}_\alpha$ and $\mathbf{S}_\alpha$ are equivalent, thanks to the functor $\Phi$ from the first to the second:

$$\Phi_{\mathrm{Obj}}(x) = e^{2i\pi x} \quad \text{and} \quad \Phi_{\mathrm{Mor}}(x, \mathbf{t}_{n+\alpha m}) = (e^{2i\pi x}, e^{2i\pi\alpha m}).$$

Moreover, they are also MRW-equivalent, by considering the set of germs of local diffeomorphisms $x \mapsto e^{2i\pi x}$, everywhere from $\mathbf{R}$ to $S^1$. Therefore, their associated $\mathbf{C}^*$-algebras are Morita equivalent. The algebra associated with $\mathbf{S}_\alpha$ has been computed numerous times and it is called *irrational rotation algebra* [MR81]. It is the universal $\mathbf{C}^*$-algebra generated by two unitary elements U and V, satisfying the relation $VU = e^{2i\pi\alpha}UV$.

REMARK 1. Thanks to the theorem (§10), and because two irrational tori $T_\alpha$ and $T_\beta$ are diffeomorphic if and only if $\alpha$ and $\beta$ are conjugate modulo $\mathrm{GL}(2, \mathbf{Z})$ [PDPI83], we get the corollary that, if $\alpha$ and $\beta$ are conjugate modulo $\mathrm{GL}(2, \mathbf{Z})$, then $\mathfrak{A}_\alpha$ and $\mathfrak{A}_\beta$ are Morita equivalent. Which is the direct sense of Rieffel's theorem [MR81, Thm 4].

REMARK 2. The converse of Rieffel's theorem is a different matter altogether. We should recover a diffeological groupoid $\mathbf{G}_\alpha$ from the algebra $\mathfrak{A}_\alpha$. Then, the space of transitive components would be the required quasifold, as stated by the proposition in (§6). In the case of the irrational torus, it is not very difficult. The spectrum of the unitary operator V is the circle $S^1$ and the adjoint action by the operator U gives $UVU^{-1} = e^{2i\pi\alpha}V$, which translates on the spectrum by the irrational rotation of angle $\alpha$. In that way, we recover the groupoid of the irrational rotations on the circle, which gives $T_\alpha$ as quasifold.

REMARK 3. Of course, the situation of the irrational torus is simple and we do not exactly know how it can be reproduced for an arbitrary quasifold. However, this certainly is the way to follow to recover the quasifold from its algebra: find the groupoid made with the Morita invariant of the algebra, which will give the space of transitivity components as the requested quasifold.

12. THE EXAMPLE OF $\mathbf{R}/\mathbf{Q}$. — The diffeological space $\mathbf{R}/\mathbf{Q}$ is a legitimate quasifold. This is a simple example with a groupoid $\mathbf{G}$ given by

$$\mathrm{Obj}(\mathbf{G}) = \mathbf{R} \quad \text{and} \quad \mathrm{Mor}(\mathbf{G}) = \{(x, \mathbf{t}_r) \mid x \in \mathbf{R} \text{ and } r \in \mathbf{Q}\}.$$

The algebra that is associated with $\mathbf{G}$ is the set $\mathfrak{A}$ of complex compact supported functions on $\mathrm{Mor}(\mathbf{G})$. Let us identify $\mathscr{C}^0(\mathrm{Mor}(\mathbf{G}), \mathbf{C})$ with $\mathrm{Maps}(\mathbf{Q}, \mathscr{C}^0(\mathbf{R}, \mathbf{C}))$ by

$$f = (f_r)_{r \in \mathbf{Q}} \quad \text{with} \quad f_r(x) = f(x, \mathbf{t}_r), \quad \text{and let} \quad \mathrm{Supp}(f) = \{r \mid f_r \neq 0\}.$$

Then,

$$\mathfrak{A} = \{f \in \mathrm{Maps}(\mathbf{Q}, \mathscr{C}^0(\mathbf{R}, \mathbf{C})) \mid \#\mathrm{Supp}(f) < \infty\}.$$

The convolution product and the algebra conjugation are, thus, given by:

$$(f * g)_r(x) = \sum_s f_{r-s}(x+s) g_s(x), \quad \text{and} \quad f_r^*(x) = f_{-r}(x+r)^*.$$

Now, the quotient $\mathbf{R}/\mathbf{Q}$ is also diffeomorphic to the $\mathbf{Q}$-circle

$$S_{\mathbf{Q}} = S^1 / \mathscr{U}_{\mathbf{Q}}, \quad \text{where} \quad \mathscr{U}_{\mathbf{Q}} = \{e^{2i\pi r}\}_{r \in \mathbf{Q}}$$

is the subgroup of *rational roots of unity*. As a diffeological subgroup of $S^1$, $\mathscr{U}_{\mathbf{Q}}$ is discrete. The groupoid $\mathbf{S}$ of the action of $\mathscr{U}_{\mathbf{Q}}$ on $S^1$ is given by:

$$\mathrm{Obj}(\mathbf{S}) = S^1 \quad \text{and} \quad \mathrm{Mor}(\mathbf{S}) = \left\{ (z, \tau) \,\Big|\, z \in S^1 \text{ and } \tau \in \mathscr{U}_{\mathbf{Q}} \right\}.$$

The exponential $x \mapsto z = e^{2i\pi x}$ realizes a MRW-equivalence between the two groupoids $\mathbf{G}$ and $\mathbf{S}$. Their associated algebras are Morita-equivalent. The algebra $\mathfrak{S}$ associated with $\mathbf{S}$ is made of families of continuous complex functions indexed by rational roots of unity, in the same way as before:

$$\mathfrak{S} = \left\{ (f_\tau)_{\tau \in \mathscr{U}_{\mathbf{Q}}} \,\big|\, f_\tau \in \mathscr{C}^0(S^1, \mathbf{C}) \text{ and } \#\mathrm{Supp}(f) < \infty \right\}.$$

The convolution product and the algebra conjugation are, then, given by:

$$(f * g)_\tau(z) = \sum_\sigma f_{\bar\sigma \tau}(\sigma z) g_\sigma(z) \quad \text{and} \quad f_\tau^*(z) = f_{\bar\tau}(\tau z)^*,$$

where $\bar\tau = 1/\tau = \tau^*$, the complex conjugate.

Now, consider $f$ and let $\mathscr{U}_p$ be the subgroup in $\mathscr{U}_{\mathbf{Q}}$ generated by $\mathrm{Supp}(f)$; this is the group of some root of unity $\varepsilon$ of some order $p$. Let $\mathrm{M}_p(\mathbf{C})$ be the space of $p \times p$ matrices with complex coefficients. Define $f \mapsto \mathrm{M}$, from $\mathfrak{S}$ to $\mathrm{M}_p(\mathbf{C}) \otimes \mathscr{C}^0(S^1, \mathbf{C})$, by

$$\mathrm{M}(z)_\tau^\sigma = f_{\bar\sigma \tau}(\sigma z), \quad \text{for all } z \in S^1 \text{ and } \sigma, \tau \in \mathscr{U}_p.$$

That gives a representation of $\mathfrak{S}$ in the tensor product of the space of endomorphisms of the infinite-dimensional $\mathbf{C}$-vector space $\mathrm{Maps}(\mathscr{U}_{\mathbf{Q}}, \mathbf{C})$ by $\mathscr{C}^0(S^1, \mathbf{C})$, with finite support.

P.I-Z — Einstein Institute of Mathematics, The Hebrew University of Jerusalem, Campus Givat Ram, 9190401 Israel.

*E-mail address*: piz@math.huji.ac.il

*URL*: http://math.huji.ac.il/~piz

E.P. — Dipartimento di Matematica e Informatica "U. Dini", Università di Firenze, Viale Morgagni 67/A, 50134 Firenze, Italy.

*E-mail address*: elisa.prato@unifi.it